\documentclass[11pt]{article}
\usepackage{pdfsync}

\usepackage[affil-it]{authblk}
\usepackage{xcolor}
\usepackage{mathrsfs}
\usepackage{amssymb,amsmath,latexsym}
\usepackage{amsfonts}
\usepackage{graphicx}
\usepackage{epstopdf}
\usepackage{caption}
\usepackage{floatrow}
\usepackage{graphicx}
\usepackage{caption}

\usepackage{atbegshi}% http://ctan.org/pkg/atbegshi
\def\I3x3{I{\ensuremath{_{3\times3}}}}

\title{On the integration of stochastic rigid body  with geometric numerical schemes}
\author{ Nataliya Ramzina\footnote{Department of Mathematical Sciences, Norwegian University of Science and Technology, 7491 Trondheim,nataliya.ramzina@math.ntnu.no}
}
\date{}

\begin{document}

%\title{
%On the integration of stochastic rigid body  with geometric numerical schemes
%}

%\author{ Nataliya Ramzina ^\star
%}
%\thanks{email: \texttt{ramzina@math.ntnu.no}}
%\email{nataliya.ramzina@math.ntnu.no}
%\affil{Department of Mathematical Sciences Norwegian University of Science and Technology, 7491 Trondheim}
%\address{$^{(1)}$Department of Mathematical Sciences Norwegian University of Science and Technology, 7491 Trondheim, nataliya.ramzina@math.ntnu.no \\
%}
%\date{\today}
%\date{March, 2014}

\maketitle

\abstract{
In the present paper we investigate the performance of explicit splitting schemes and related techniques applied to a
rigid body model subject to a stochastic torque and random perturbations in the inertia tensor. Results are discussed  and compared with traditional solvers for such model.}

%\begin{keywords}
%stochastic rigid body, splitting method, variation of parameters
%\end{keywords}

\section{Introduction}
% general staff
Stochastic differential equations (SDEs) are used in many fields, such as stock market, financial mathematics, stochastic controls,
biological science, chemical reactive kinetics and hydrology, as described in \cite{platen}, \cite{oksendal}. Thus, it is of importance to study the solution of SDEs.
However, finding the analytic solutions is often very difficult or impossible. Numerical methods, on the other hand, help to find approximate solutions of SDEs quite fast, and thus, development of such methods looks promising. However, it is of importance to study the solution of SDEs obtained numerically, because of results obtained by different methods may deviate significantly in terms of accuracy and computational resources required. 
%Practical applications are usually consider a rigid bodies subject to stochastic forces and moves. 
  
%In this paper we consider a stochastic rigid body model. Imagine the body is subject to small random perturbations such that each of them causes a displacement of the body within time. 

Many of the physical systems could be refined in terms of generalized stochastic rigid body model. For example, this model might be rigid body subject to small perturbation force randomly dependent on time, or the body with randomly perturbed inertia. Examples of such an approach could be found in polymer simulations when the polymer chain is divided into several interacting monomers treated as rigid bodies during the modelling. The stochastic torques arise when the polymer is immersed in a solvent: the kinematics of each of the bodies in the model is influenced by stochastic loads raised from thermal fluctuations \cite{polymer_dynamics}.

%Such effects might be a result for example from a solvent
%when the polymer, immersed in it, modelled as interacting rigid bodies. Kinematics of such each body is %influenced by  stochastic  loads (\cite{polymer_dynamics}).
%A relevant example is the one of a 
 Another example is considered in the work by Arnold et al. \cite{Random_Dynamical} where authors concentrate on a ship partly immersed in water. The authors propose a stochastic model which includes an impact of severe weather conditions such as seaway and wind on a roll motion of the ship. 

Another example arises in the field of animation \cite{xie}, where a light body flying through air is simulated by solving generalized Kirchhoff equations. 
The challenge is to determine resulted forces and torques due to the turbulence generated at the body surface (vortical loads). In order to do that effectively the authors suggest a coupled approach: solve the free rigid body equations with no vortical loads and  then generalized  Langevin equations used to represent the characteristics of the surrounding flow. 

Each of the above examples is treated with specifically tailored numerical schemes.

% our work

In this work we focus on the evolution of  the angular momentum of a rigid body subject to a stochastic torque. This simple test case is used here for testing our numerical schemes. We also study a rigid body model with randomly perturbed inertia tensor. This is of interest in the case of multibody system dynamics with uncertain rigid bodies as studied in \cite{random_inertia}.

 We consider two techniques:  a splitting method and a method based on the variation of constants formula. These are well known techniques for deterministic problems which we here want to test in a stochastic context. We believe these techniques have some potential also for more complex problems \cite{random_inertia}, \cite{xie}, \cite{Random_Dynamical} then the ones considered here. 
%We also study a rigid body model with randomly perturbed inertia tensor. This is of interest in the case of multibody system dynamics with uncertain rigid bodies as studied in \cite{random_inertia}.

 We start with the formulation of the mathematical model followed by the description of the numerical schemes. 
Then we  apply  a variety of numerical schemes and investigate them in terms of weak errors. We compare the  methods in terms of  CPU time and relative cost. Further we investigate the geometrical properties of the numerical methods applied to the model with randomly perturbed inertia.

%%%%%%%%%%%%%%%%%%%%%%%%%%%%%%%%%%%%%%%%%%%%%%%%%%%%%%%%%%%%%%%%%%%%%%%%%%%%%%%%%%%%%%%%%%%%%%%%%%%%%%%%%%%%%%%%%%%%%%%%%%%%%%
\section{Rigid body under random perturbation}
\subsection{Mathematical model}
We study the following stochastic differential equation with one-dimensional Brownian motion

%\begin{equation}\label{eq:full_sde}
%\left\{\begin{array}{lcl}d \boldsymbol{m} &=& (\boldsymbol{m} \times T^{-1}\boldsymbol{m})dt+a\,\boldsymbol{m}\,dW_{t} ,\qquad \boldsymbol{m}(0)=\boldsymbol{m}_0\\
%d{\mathbbm{q}} &=& \frac{1}{2}\mathbbm{q}\,\Omega dt , \qquad \Omega=(0,\boldsymbol{\omega}) 
%\end{array} \right. ,
%\end{equation}
\begin{equation}\label{eq:full_sde}
d \boldsymbol{m} = (\boldsymbol{m} \times T^{-1}\boldsymbol{m})dt+a\,\boldsymbol{m}\,dW_{t} ,\qquad \boldsymbol{m}(0)=\boldsymbol{m}_0
\end{equation}
where $\boldsymbol{m} \in  \mathbb{R}^3 $ is the angular momentum of the body, $T \in  \mathbb{R}^{3\times 3} $ is a diagonal inertia matrix, $\boldsymbol{\omega}=T^{-1}\boldsymbol{m}\in  \mathbb{R}^3  $ is angular velocity,
$W_t$ is a standard Brownian motion, $a$ is a parameter (noise). We assume that the three moments of inertia are pairwise distinct and we order them in ascending order.

% $\mathbbm{q} \in  \mathbb{R}^4$ is a quaternion representing an attitude of the body, 
% defined on a filtered probability space $(\Omega,\, \mathscr{F},\, \lbrace\mathscr{F}_t\rbrace_{t\in [0,T]}, \boldsymbol{P})$, $a$ is a white one-dimensional noise. Let us give a definition of filtered probability space. Let $(\Omega, \, \mathscr{F}, \, \boldsymbol{P})$ be a probability space. Here $\Omega$ is the sample space, $\mathscr{F}$ is a $\sigma$ - algebra on $\Omega$, and $\boldsymbol{P}$ a probability measure on $\Omega$. A filtration $\lbrace\mathscr{F}_t\rbrace_{t\in [0,T]}}$ in $(\Omega, \, \mathscr{F}, \, \boldsymbol{P})$ is an increasing family of $\sigma$-algebras such that 
%$$\mathscr{F}_0\subset\mathscr{F}_t\subset\dots\mathscr{F}_T\,=\mathscr{F}.$$
%Then $(\Omega,\, \mathscr{F},\, \lbrace\mathscr{F}_t\rbrace_{t\in [0,T]}, \boldsymbol{P})$ is called a \emph{filtered probability space}.

%%%%%%%%%%%%%%%%%%%%%%%%%%%%%%%%%%%%%%%%%%%%%%%%%%%%%%%%%%%%%%%%%%%%%%%%%%%%%%%%%%%%%%%
\subsection{Numerical methods}
\subsubsection{Splitting scheme}

Since \eqref{eq:full_sde}  is a sum of two exactly solvable terms  and splitting it into deterministic and stochastic parts gives us exact solution

%\begin{equation}\label{eq:flows1}
%S_1 = \left\{\begin{array}{lcl}d \boldsymbol{m} &=& (\boldsymbol{m} \times T^{-1}\boldsymbol{m})dt, \\
%d\mathbbm{q} &=& \frac{1}{2}\mathbbm{q}\,\Omega \,dt 
%\end{array} \right. ,
%\end{equation}
%and
%\begin{equation}\label{eq:flows2}
%S_2 = \left\{\begin{array}{lcl}d \boldsymbol{m} &=&a\,\boldsymbol{m}\,dW_{t} \\
%d\mathbbm{q} &=& 0 \\
%\end{array} \right. 
%\end{equation}
\begin{equation}\label{eq:flows1}
S_1\,:\qquad d \boldsymbol{m} = (\boldsymbol{m} \times T^{-1}\boldsymbol{m})dt
\end{equation}
\begin{equation}\label{eq:flows2}
S_2\,:\qquad d \boldsymbol{m} =a\,\boldsymbol{m}\,dW_{t}
\end{equation}
then we compose the flows of systems  \eqref{eq:flows1} and \eqref{eq:flows2} into the first order splitting scheme
 \begin{equation}
\label{trotter}
(\boldsymbol{m}, q)^{(j+1)} = \phi_{h}^{[S_2]}\circ \phi_{h}^{[S_1]} \,
((\boldsymbol{m}, q)^{(j)}).
\end{equation}

The first equations \eqref{eq:flows1} are a free rigid body motion equations and could be integrated exactly. An explicit solution can be
obtained by using Jacobi elliptic functions and has been implemented to machine accuracy following the techniques presented in  \cite{celledoni-2008}. The solution of the second set of equations \eqref{eq:flows2} is

\begin{equation}\label{flow2:exact}
\boldsymbol{m}(t) =\exp\left(
a\,W(t)-\frac{1}{2}\,a^2t
%\begin{array}{ccc}
%0 & 0 & 0 \\
%0 & 0 &a\,dW(t)\\
%0 & -a\,dW(t) & 0
%\end{array}
\right) \boldsymbol{m}_0 \end{equation}

where $W(t)$ is independent random variable of the form $\sqrt{t}N(0,1)$.

%%%%%%%%%%%%%%%%%%%%%%%%%%%%%%%%%%%%%%%%%%%%%%%%%%%%%%%%%%%%%%%%%%%%%%%%%%%%%%%%%%%%%
\subsubsection{The Variation of  Constants formula}
We recall the variation of constants formula. 
Let us consider the following system of differential equations
\begin{equation}\label{v}
\begin{aligned}
\dot{\boldsymbol{z}}&=f(\boldsymbol{z})+g(\boldsymbol{z}),&\boldsymbol{z}(0)&=y_0,
\\
\dot{\boldsymbol{y}}&=f(\boldsymbol{y}), & \boldsymbol{y}(0)&=y_0\\
\end{aligned}
\end{equation}
and suppose that $\frac{\partial f}{\partial \boldsymbol{y}}$ exists and is continuous. Then the solution is given by the non-linear variation of constants formula  \cite{hairer} (p.96, Chapter 1).
\begin{equation}\label{voc}
\boldsymbol{z}(t)=\boldsymbol{y}(t)+\int_{0}^{t}\Phi(t,\,s,\,\boldsymbol{z}(s))\cdot g(\boldsymbol{z}(s))ds,
\end{equation}
with $\Phi$ is the solution of the equation 

\begin{equation}\label{resolvent}\begin{array}{lcl}
\Phi^\prime=\frac{\partial f}{\partial \boldsymbol{y}} (\boldsymbol{y})\cdot\Phi,\\
\Phi(0)=\I3x3\\
\end{array}
\end{equation}
solved together with  $\dot{\boldsymbol{y}}=f(\boldsymbol{y}) $.
For our problem in order to integrate \eqref{eq:full_sde} with the variation of constants formulae we set
\begin{equation}\label{part}
\dot{\boldsymbol{y}}=f(\boldsymbol{y})=\boldsymbol{y}\times T^{-1}\boldsymbol{y},
\end{equation} and 
\begin{equation}
g\,dt=a\,\boldsymbol{z}dW_t.
\end{equation}
Let us denote with 
$$A(t,\boldsymbol{y})=\frac{\partial f}{\partial \boldsymbol{y}} (\boldsymbol{y})$$ then we obtain
$$A(t,\boldsymbol{z})=\left(\begin{array}{ccc}
0&\boldsymbol{z}(3)(\frac{T_2-T_3}{T_2T_3})&\boldsymbol{z}(2)(\frac{T_2-T_3}{T_2T_3})\\[6pt]
\boldsymbol{z}(3)(\frac{T_3-T_1}{T_1T_3})&0&\boldsymbol{z}(1)(\frac{T_3-T_1}{T_1T_3})\\[6pt]
\boldsymbol{z}(2)(\frac{T_1-T_2}{T_2T_1})&\boldsymbol{z}(1)(\frac{T_1-T_2}{T_2T_1})&0
\end{array}\right).$$ 
Then we use \eqref{voc} to construct our numerical approximation 
$\boldsymbol{z}_{n+1}=\boldsymbol{y}_{n+1}+a\, dW_{n}\, \Phi_{n+1} \boldsymbol{z}_{n}$, where $\boldsymbol{y}_{n+1}$ is a solution of \eqref{part}, $\Phi_{n+1}$ is a solution of \eqref{resolvent}, $dW_n$ is a Brownian increment. Equation \eqref{part} is a free rigid body problem and we solve it exactly with Jacobi elliptic functions \cite{celledoni-2008}, obtaining $\boldsymbol{y}_1=\boldsymbol{y}(t_1),\,\dots,\, \boldsymbol{y}_n=\boldsymbol{y}(t_n)$. 
We integrate \eqref{resolvent} with a Magnus method of order two \cite{hairerg} 
and get
\begin{equation}\label{magnus}
\Phi_{n+1}=\exp(hA(t_n+h/2, {\boldsymbol{z}_{n+1/2}}))\Phi_n,
\end{equation}
%where ${\boldsymbol{z}_{n+1/2}}$ is an approximation of ${\boldsymbol{z}(t)$ for $t=\frac{t_n+ t_{n+1}}{2}$. The numerical solution ${\boldsymbol{z}_{n+1}\approx{\boldsymbol{z}(t_{n+1})$ is 
%$${\boldsymbol{z}_{n+1}={\boldsymbol{y}_{n+1}+a\,dW_{n}\,\Phi_{n+1}\boldsymbol{z}_{n},$$ 
where 
$$\boldsymbol{z}_{n+1/2}=\boldsymbol{y}_{n+1/2}+a\,dW_{n}\,\Phi_{n+1/2}\boldsymbol{z}_{n},\quad \Phi_{n+1/2}=\exp ( h A( t_n+h/2, \boldsymbol{y}_{n+1/2} ))\Phi_n.$$
%We proceed as following hence: \begin{itemize}
%\item solve \eqref{part} for $t=t_n+h/2$
%\item evaluate \eqref{magnus} at this solution and freeze it
%\item solve \eqref{voc} to find $\boldsymbol{z}(t_n+h/2)$;
%\item evaluate \eqref{magnus} along  the trajectory $\boldsymbol{z}(t_n+h/2)$ and freeze it;
%\item solve \eqref{part} for $t=t_n+h$;
%\item finally get the solution $\boldsymbol{z}(t_n+h)$ from \eqref {voc}
%\end{itemize}

%Then insert output from \eqref{part} at time $t_{n+1}$ and solution of \eqref{resolvent} into
%\begin{equation}
%\boldsymbol{z}(t_{n+1})=\boldsymbol{y}(t_{n+1})+h\sum_{j=0}^TdW_j\Phi(t_{n+1}-s_j, \, \boldsymbol{z}(s_j))g(\boldsymbol{z}_j)
%\end{equation}

%%%%%%%%%%%%%%%%%%%%%%%%%%%%%%%%%%%%%%%%%%%%%%%%%%%%%%%%%%%%%%%%%%%%%%%%%%%%%%%%%%%
\section{Random perturbation of inertia tensor}
Consider a free rigid body model subject to white noise in inertia 
\begin{equation}\label{eq:full_sde2}
\dot{\boldsymbol{m}} = \boldsymbol{m} \times T^{-1}\boldsymbol{m} ,\qquad \boldsymbol{m}(0)=\boldsymbol{m}_0
\end{equation}
%\begin{equation}\label{eq:full_sde2}
%\left\{\begin{array}{lcl} \dot{\boldsymbol{m}} &=& \boldsymbol{m} \times T^{-1}\boldsymbol{m} ,\qquad \boldsymbol{m}(0)=\boldsymbol{m}_0\\
%\dot{\mathbbm{q}} &=& \frac{1}{2}\mathbbm{q}\,\Omega , \qquad \Omega=(0,\boldsymbol{\omega}) 
%\end{array} \right. ,
%\end{equation}
with $T=T_d+T_s$, $T_d$ is a deterministic inertia tensor, $T_s=\epsilon\,W$ with a small parameter $\epsilon$ and an arbitrary martingale 
$$W=\left(\begin{array}{ccc}
\xi_{11}&\xi_{12}&\xi_{13}\\
\xi_{21}&\xi_{22}&\xi_{23}\\
\xi_{31}&\xi_{32}&\xi_{33}
\end{array}\right),\quad \xi_{ij}\in \mathbb{R}$$

We then need to calculate an inverse of the inertia matrix 
\begin{equation*}
(T_d+T_s)^{-1}=T_d^{-1}(\I3x3+T_sT_d^{-1})^{-1}.
\end{equation*}
%$\Vert T_sT_d^{-1}\Vert<1$
Assuming $\Vert T_sT_d^{-1}\Vert<1$ with $\Vert \bullet\Vert$ any matrix norm, we can apply Neumann series (for which convergence guaranteed) and we get
\begin{equation*}
T_d^{-1}\sum_{k=0}^{\infty}(T_sT_d^{-1})^k=T_d^{-1}(\I3x3+\sum_{k=1}^{\infty}(T_sT_d^{-1})^k)=T_d^{-1}+\sum_{k=1}^{\infty}(T_sT_d^{-1})^k.
\end{equation*} 
Assume $\tilde{T_d}=T_d^{-1}$ and $\tilde{T_s}=\sum_{k=1}^{\infty}(T_sT_d^{-1})^k,$ such that 
\begin{equation}\label{eq:ri_sde}
\left\{\begin{array}{lcl} \dot{\boldsymbol{m}} &=& \boldsymbol{m} \times \tilde{T_{d}}\boldsymbol{m} +\boldsymbol{m} \times \tilde{T_{s}}\boldsymbol{m} ,\\
\dot{q} &=& \frac{1}{2}q\,\Omega 
\end{array} \right. 
\end{equation}

%%%%%%%%%%%%%%%%%%%%%%%%%%%%%%%%%%%%%%%%%%%%%%%%%%%%%%%%%%%%%%%%%%%%%%%%%%%%%%%%%%%%%
\subsection{Numerical solution}
In the case of random perturbation of inertia tensor we split \eqref{eq:ri_sde} into stochastic and deterministic parts 
%\begin{equation}\label{eq:flows11}
%S_1 = \left\{\begin{array}{lcl}\dot{\boldsymbol{m}} &=&  \boldsymbol{m} \times \tilde{T_{d}}\boldsymbol{m} , \\
%\dot{\mathbbm{q}} &=& \frac{1}{2}\mathbbm{q}\,\Omega 
%\end{array} \right. ,
%\end{equation}
\begin{equation}\label{eq:flows11}
S_1\,:\qquad \dot {\boldsymbol{m}}= \boldsymbol{m} \times \tilde{T_{d}}\boldsymbol{m}
\end{equation}
and
\begin{equation}\label{eq:flows22}
S_2\,:\qquad \dot {\boldsymbol{m}} = \boldsymbol{m} \times \tilde{T_{s}}\boldsymbol{m}
\end{equation}
%\begin{equation}\label{eq:flows22}
%S_2 = \left\{\begin{array}{lcl}\dot {\boldsymbol{m}} &=&\boldsymbol{m} \times \tilde{T_{s}}\boldsymbol{m}, \\
%\dot{\mathbbm{q}} &=& 0 \\
%\end{array} \right. 
%\end{equation}
We then combine the two flows into Lie-Trotter scheme.

We shall also use solution from (fully) implicit midpoint rule given by the general scheme in stochastic setting \cite{milstein}

\begin{equation}
X_{n+1}=X_n+h\,f\left(\frac{X_n+X_{n+1}}{2}\right),
\label{eq:imr}
\end{equation}

where $X_n$ is a numerical approximation of $\dot{X}=f(X)$.
%%%%%%%%%%%%%%%%%%%%%%%%%%%%%%%%%%%%%%%%%%%%%%%%%%%%%%%%%%%%%%%%%%%%%%%%%%%%%%%%%%%
 
\section{Results and discussion}
\begin{figure}[ht!]
\includegraphics[scale=0.6]{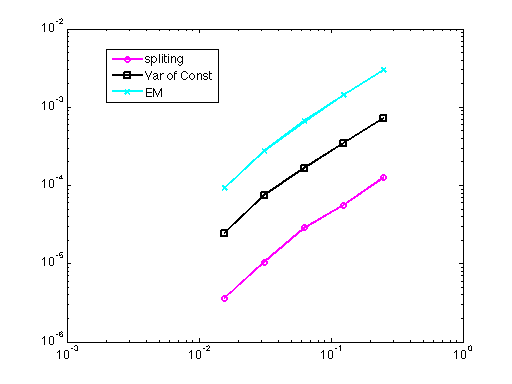}
\caption{Error $\vert \mathbb{E}(\boldsymbol{m}(T))-\mathbb{E}(\boldsymbol{m}_{ref})\vert$  where $\boldsymbol{m}$ a  solution of \eqref{eq:full_sde} with \emph{Lie-Trotter splitting} the Euler-Maruyma and Variation of Constants.}
\label{spl_EM}
\end{figure}

Firstly we show the weak convergence of the splitting method, the Euler-Maruyma (EM)  and the Variation-of-Constants (VoC) methods applied to problem \eqref{eq:full_sde}.
Figure~\ref{spl_EM}  shows the absolute error between mean values over sampled paths $\vert \mathbb{E}(\boldsymbol{m}(T))-\mathbb{E}(\boldsymbol{m}_{ref})\vert$ in the angular momentum versus time integration step. Initial values $\boldsymbol{m}_0$ = [0.4165;0.9072;0.0577], T= diag([0.9144;1.098;1.66]), $a=0.1$,  integration time t = 1. The number of integration steps have been set ranging from $2^1$ up to $2^{9}$ and equally spaced in logarithmic scale. S = 1000 paths are sampled. For a reference solution we take the approximated one with step-size $2^{-8}$. 

\begin{figure}[ht!]
\includegraphics[scale=0.6]{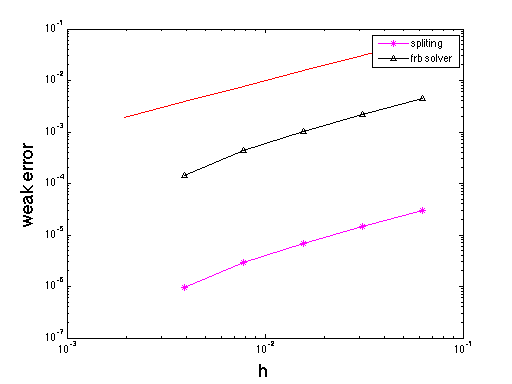}
\caption{Error $\vert \mathbb{E}(\boldsymbol{m}(T))-\mathbb{E}(\boldsymbol{m}_{ref})\vert$ with  $\boldsymbol{m}$ a solution of \eqref{eq:full_sde2} with \eqref{trotter} (magenta) and FRB integrator (black). Reference red line has slope one.}
\label{random_inertia}
\end{figure}

On Figure~\ref{random_inertia} one can see the weak convergence of splitting solver \eqref{trotter} applied to \eqref{eq:flows11}, \eqref{eq:flows22} and the free rigid body (FRB) integrator
\cite{celledoni-2008}. 
We note the error in reference solution significantly depends on number of paths sampled, see Table~\ref{M_TABLE}.

\begin{center}
\begin {table}[h!]
\centering
\caption {Absolute error in reference solution vs number of Brownian paths sampled}\label{M_TABLE}
%\begin{tabular}{ | p{6cm} | C{2.5cm} | C{2.5cm} | }
\begin{tabular}{ | c | c | c | c | }
\hline
S       & Error     \\ \hline
100    & 0.0016     \\ \hline
1000  & 4.38e-04   \\ \hline
10000  & 6.38e-05   \\ \hline
100000 & 7.71e-05     \\ \hline

\end{tabular}
\end{table}
\end{center}

\begin{figure}[h]
      
           \includegraphics[width=0.8\textwidth]{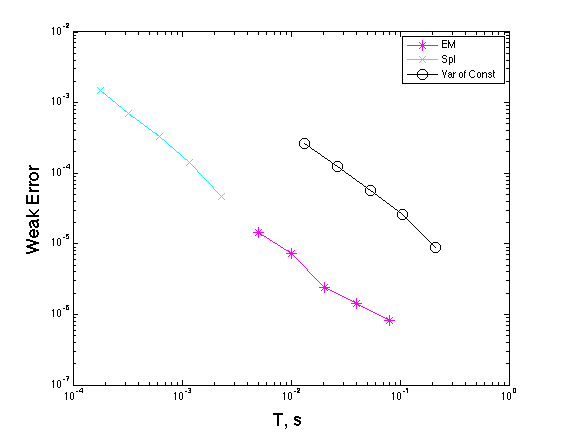}
           \caption{Average log of error versus cpu time (1000 Brownian paths).}
       \label{fig:cpu}

\end{figure}

From the above Figure~\ref{spl_EM} it is easy to see that the EM scheme, the splitting scheme and the VoC have the order one weak convergence, but the splitting scheme obtains better approximate solutions for the problem \eqref{eq:full_sde}.
The evaluation of CPU time used for the calculations have been performed directly in MATLAB\textregistered \,by measuring the time interval between the initiation of the algorithm and completion of the calculation. The performance of the algorithm is well represented by the plot of the weak error in  $\boldsymbol{m}$  against the CPU time (Figure~\ref{fig:cpu}). The amount of calculations is proportional to the number of integration steps hence we expect the EM to be the fastest when no high accuracy is required. 

The relative cost of the method can be computed as a ratio between given method cost and the minimal cost of method.
%\begin{equation}
%\frac{Ac(X)}%{min_{all\thinspace methods} aver.\thinspace cost\thinspace of\thinspace method\thinspace X}
%{min_{S}Ac(X)}
%\label{eq: costi}
%\end{equation}
%where Ac is an average cost of method  X, S is a set of the studied methods.
The cost of Euler-Maruyma method is  minimal compare to the Variation of Constants and splitting scheme. We see that the relative cost of the splitting scheme and variation of constants is decreasing for large step-sizes. This makes them attractive when large step-sizes are required in simulations.

\begin{figure}[h]

        \includegraphics[width=0.8\textwidth]{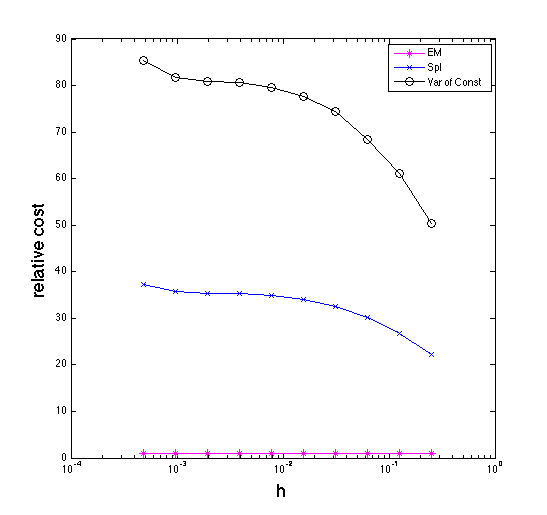}
        
        \caption{Relative cost (with respect to the cheapest method) versus step-size}
           \label{fig:cost}

              \end{figure}

\begin{figure}[ht]
\includegraphics[width=1\textwidth]{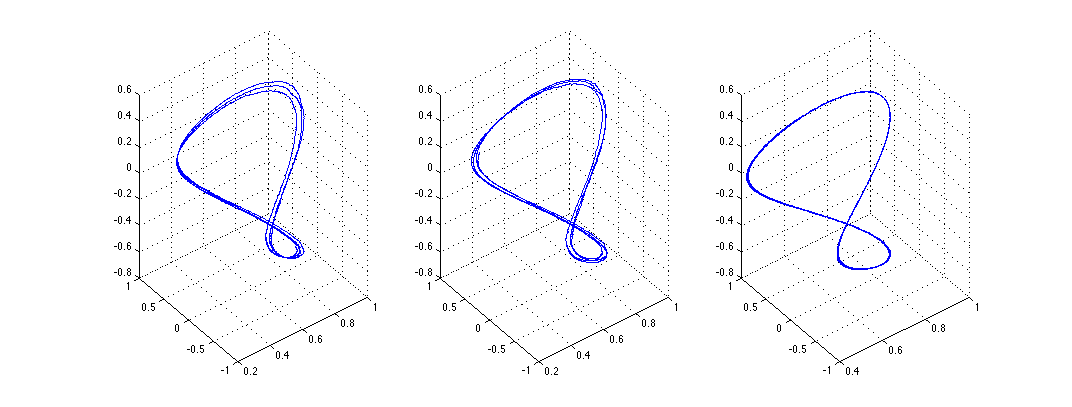}
\caption{The vector $\boldsymbol{m}_i \approx \boldsymbol{m}(ih)$ for the perturbed rigid body motion over  time
 T=10000 with h=0.1, noise parameter a=0.001. Left splitting method; middle EM method; right variation of constants  \label{fig:smallest_noise} }
\end{figure}

\begin{figure}[ht]
\includegraphics[width=1\textwidth]{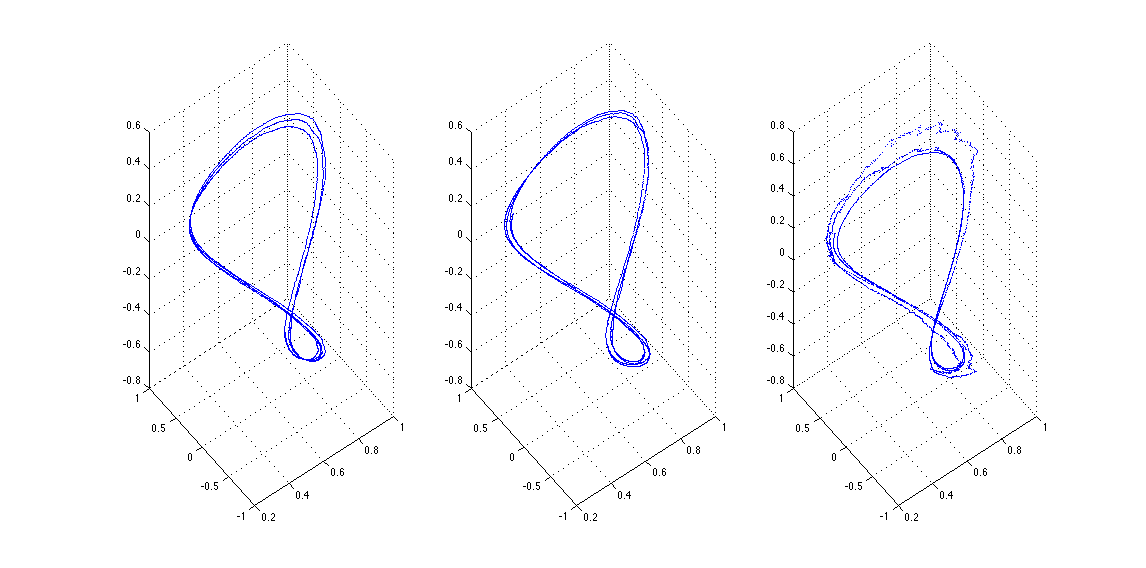}
\caption{The vector $\boldsymbol{m}_i \approx \boldsymbol{m}(ih)$ for the perturbed rigid body motion over  time
 T=10000 with h=0.1, noise parameter a=0.005. Left splitting method; middle EM method; right variation of constants  \label{fig:small_noise}}
\end{figure}

We next consider a perturbed rigid body where the perturbation is given by a random torque. A similar problem has been considered in the deterministic case for example in \cite{C_Z}.  The aim of the experiments is to see how the trajectories of the momentum are affected by the size of the noise. We see that the numerical trajectory maintains the same character of the free rigid body motion when the noise parameter $a$ is small enough see Figure~\ref{fig:smallest_noise}.
This behaviour is due to usage of free rigid solver based on Jacobi elliptic functions \cite{C_Z}, \cite{celledoni-2008}, \cite{niclas}.  However when the noise term starts increasing the trajectory of variation of constants method does not look smooth anymore, see Figure~\ref{fig:small_noise}.
%%% Choose the middle noise parameter (for which the character is still preserved) and plot figures with 
%%% that one for the variation of constant and EM

%\begin{figure}
%\label{fig:big_torque}
%\includegraphics[width=1\textwidth]{Figures/big_noise.eps}}\
%\caption{The vector $\boldsymbol{m}_i \approx \boldsymbol{m}(ih)$ for the perturbed rigid body motion over  time
% T=10000 with h=0.1, a=0.1}
%\end{figure} 

\begin{figure}[h]
\centering
        \includegraphics[width=1\textwidth]{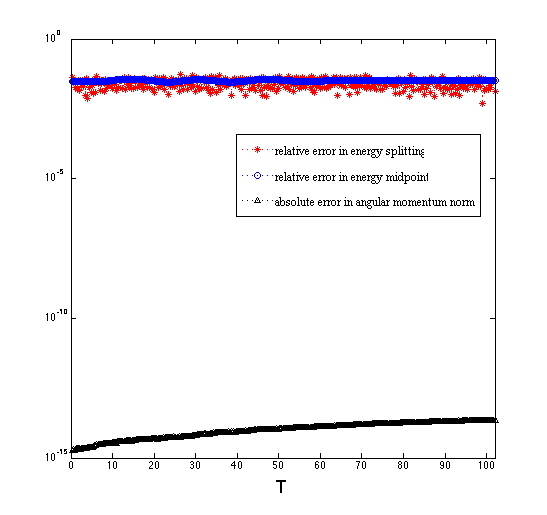}
        
        \caption{Energy error for the body with random inertia. Integration on the interval [0,100], \emph{h}=0.25, using the splitting method  and the midpoint rule. \label{fig:geomp}}
           
\end{figure}

Another issue, that might be of interest is the preservation of angular momentum in models with random inertia. Also we look at energy of the model and see to which extent it preserved. We compare our results with the implicit midpoint rule and present them in Figure~\ref{fig:geomp}. We measure the energy as 
\begin{equation}
H=\frac{1}{2} \omega^T \,T\,\omega,
\label{eq:energy}
 \end{equation} 
 with  $\boldsymbol{\omega}=T^{-1}\boldsymbol{m}\in  \mathbb{R}^3  $, $T=T_d+T_s$. The energy error is averaged over 100 Brownian samples.  

\section{Conclusion}
 In the current work we have performed numerical methods, which widely used for deterministic problems, to simple stochastic models. We see that methods work the same way as expected. We propose that these well know methods have a prospectiveness for more 
advanced stochastic models.

% In the current we have performed  stochastic modelling by combing widely used methods for NORMAL SOLVERS incorporated with STOCHASTIC PART.

% Thus the modelling becomes more general and system-independent. Therefore the advantage of this model's usage is reduced human resources required to model the given problem. 

%In the current work we applied splitting schemes to resolve the dynamics of a marine vessel. We compared the CPU time required to accomplish the calculations with the relative error of calculations. The results suggest that higher order splitting schemes achieve the aimed precision faster  compared to classical Runge-Kutta schemes. The use of splitting schemes became more popular in scientific simulation in the past decade and this study confirms there is a potential for the use of such schemes also in engineering problems.

\section{Acknowledgements}

This research was supported by a Marie Curie International Research Staff Exchange Scheme Fellowship within the 7th European Community Framework Programme. The authors would like to acknowledge the support from the GeNuIn Applications project funded by the Research Council of Norway, and part of the work was carried out while the authors were visiting Massey University, Palmerston North, New Zealand and La Trobe University, Melbourne, Australia.


\begin{thebibliography}{99}
%
\bibitem{blanes-2002}
S. ~Blanes and P.C.~ Moan.
Practical symplectic partitioned Runge-Kutta and Runge-Kutta-Nystr\"{o}m methods
\textit{Journal of Computational and Applied Mathematics},
142:313--330, 2002.
%
\bibitem{Random_Dynamical}
L.~Arnold, I.~Chueshov, and G.~Ochs.
Random Dynamical Systems Methods in Ship Stability: A Case Study  
{\em Journal of Physics A}, 39:5463--5478, 2006.
 
 %
 \bibitem{random_inertia}
A.~Batou, C.~Soize
\newblock {\em Multibody system dynamics with uncertain rigid bodies}.
\newblock  Proceedings of the 8th International Conference on Structural Dynamics, EURODYN 2011, 2620--2625. 

 
%
\bibitem{celledoni-2008}
E. ~Celledoni, F. ~Fass\'{o}, N.~ S\"{a}fstr\"{o}m and A.~ Zanna.
The exact computation of the free rigid body motion and its use in splitting methods
\textit{SIAM Journal on Scientific Computing},
30:4:2084--2112, 2008.
%
\bibitem{niclas}
E. ~Celledoni and N.~ S\"{a}fstr\"{o}m.
\newblock Efficient time-symmetric simulation of torqued rigid bodies using Jacobi elliptic functions.
\newblock Journal of Rhysics A, Vol. 39, pp. 5463--5478, 2006
%
\bibitem{C_Z}
E. ~Celledoni and A.~ Zanna.
FRB-FORTRAN routines for the exact computation of free rigid body motions
\textit{ACM ToMS}, vol 37 nr. 2 (2010).
%
\bibitem{chaturvedi-2011}
N.A.~ Chaturvedi,  A.K.~ Sanyal and N.H.~ McClamroch
Rigid-Body Attitude Control
\textit{IEEE Control Systems Magazine},
31:3:30--51, 2011.
%
\bibitem{dullweber-1997}
A. Dullweber, B. Leimkuhler, and R. McLachlan.
Symplectic splitting methods for rigid body molecular dynamics.
\textit{J. Chem. Phys},
107:5840--5851, 1997.
%
\bibitem{fossen-2002}
T.I. ~Fossen.
\textit{Marine Control Systems}.
Marine Cybernetics, 2002.
%
\bibitem{geradin-2001}
M. Geradin and A. Cardona.
\textit{Flexible Multibody Dynamics}.
Wiley and Sons Ltd., 2001.
%
\bibitem{hairerg}
E.~Hairer, C.~Lubich, and G.~Wanner.
\newblock Geometric numerical integration, volume 31 of Springer series in computational mathematics.
\newblock Springer series in computational mathematics, volume 31. Springer, 2002.

%
\bibitem{hairer}
E.~Hairer, S.P.~N\o rsett and G.~Wanner.
\newblock{\em Solving Ordinary Differential equations. Nonstiff Problems}.
\newblock Springer series in computational mathematics 8. Springer, Berlin, 2002.

%
\bibitem{milstein}
G.N.~Milstein, Yu.M.~Repin, and M.V.~Tretyakov.
\newblock Numerical Methods for Stochastic Systems Preserving Symplectic Structure.
\newblock SIAM J. Numer. Anal. Vol. 40, No. 4, pp.1583--1604.
%
\bibitem{perez-2007}
T. Perez and T.I. Fossen.
Kinematic Models for Manoeuvring and Seakeeping of Marine Vessels
\textit{Modeling, Identification and Control}, 
28:1:19--30, 2007.
%

\bibitem{platen}
E. ~Platen
\newblock{\em An introduction to numerical methods for stochastic differential equations}.
\newblock Acta Numer., 8 (1999), pp. 197--246.

%
\bibitem{polymer_dynamics}
J.~Walter, O.~Gonzalez and J.H~Maddocks.
\newblock {\em On the stochastic modelling of rigid body systems with application to polymer dynamics}.
\newblock  SIAM Multiscale Model. Simul., Vol.8, No.3, 1018--1053, 2010

%
\bibitem{xie}
H.~Xie, K.~Miyata.
\newblock {\em Stochastic modelling of immersed rigid-body dynamics}.
\newblock Proceedings of SIGGRAPH Asia 2013
%
\bibitem{oksendal}
B.~{{\O{}}}ksendal 
\newblock{\em Stochastic differential equations: an introduction with applications}.
\newblock Springer, 6th edition, 2010

\end{thebibliography}
\end{document}